\documentclass[12pt]{article}

\usepackage{latexsym}
\usepackage[utf8]{inputenc}
\usepackage{graphicx}
\usepackage{amsmath}
\usepackage{dsfont}
\usepackage{amssymb}
\usepackage{amsthm}
\usepackage{color}

\usepackage{amsfonts}
\usepackage{latexsym}
\usepackage{amssymb}

\newtheorem{theorem}{Theorem}[section]
\newtheorem{proposition}[theorem]{Proposition}
\newtheorem{corollary}[theorem]{Corollary}

\def\E{{\mathds E}}
\def\N{{\mathds N}}

\def\R{{\mathds{R}}}

\def\l{{\lambda}}

\begin{document}

\title{Best and random approximation \\ of a convex body by a polytope}
\author{J. Prochno, C. Sch\"utt, E. M. Werner}
\maketitle

\begin{abstract}
In this paper, we give an overview of some results concerning
best and random approximation of convex bodies by polytopes. We explain how both are linked and see that random approximation is almost as good as
best approximation.
\end{abstract}
\text{Primary: 52A22, 52A27, 60D05}
\text{Secondary: 52A20}\\
\text{Keywords: affine surface area, convex body, Delone triangulation constant, polytope}

\section{Introduction}
A common and important question in mathematics is whether
mathematical objects or constructions with certain
requirements or features exist. Such questions naturally appear in various areas of mathematics and theoretical computer science. For instance, in numerical linear algebra, where randomization is used for designing fast algorithms \cite{KV2017}, in algorithmic convex geometry, in relation to the complexity of volume computation of high-dimensional convex bodies \cite{AGB2015,V2010}, in geometric functional analysis, when one is interested in constructing normed spaces with certain pathological characteristics, or in graph theory, where one is interested in finding graphs with certain pre-described features. We will elaborate on the latter two and further examples below, where the corresponding references will be provided.
\par
With a view towards motivation, let us stick with the last topic in that list for a moment, graph theory. The chromatic number of a graph is the least amount of colors necessary to color the vertices such that no two adjacent vertices share the same color and the girth of a graph is the length of a shortest cycle contained in
the graph. In 1959, Erd\H{o}s proved that there are graphs whose girth and chromatic numbers are both arbitrarily large \cite{E1959}. Another classical result in graph theory is the following: there is a constant $c\in(0,\infty)$ such that for all sufficiently large $n\in\N$ there exists a graph with $n$ vertices which contains no triangle and which does not
contain a set of $c\sqrt{n} \ln n$ independent vertices; recall that a set of independent vertices in a graph is a set of vertices such that
no two vertices represent an edge of the graph. In this case, Erd\H{o}s could not give an explicit construction, but he showed that such graphs exist with high probability \cite{E1961}.
\par
Erd\H{o}s' probabilistic approach is commonly referred to as the probabilistic method today and he was arguably the most famous of its pioneers (although others before him proved theorems using this method). Also Shannon employed this method in the proof of his famous source coding theorem \cite{S1948}. The probabilistic method has become mathematical folklore way beyond graph theory to use it in order to show the existence of objects with prescribed features.
\par
Let us continue our motivation on more analytic grounds. A famous theorem of Dvoretzky \cite{Dvo1961} states that in every 
infinite dimensional Banach space there are finite dimensional subspaces of
arbitrarily large dimension that are up to a small error 
Euclidean spaces. More specifically, let $X$ and $Y$ be two normed spaces.
The Banach-Mazur distance between $X$ and $Y$ is
$$
d(X,Y):=\inf\big\{\|A\|\|A^{-1}\|\,:\, A:X\to Y\hskip 2mm\mbox{is an isomorphism}\big\},
$$
and, in case the spaces are not isomorphic, $d(X,Y):=+\infty$.
Dvoretzky's theorem says: there is a constant $c\in(0,\infty)$ such that for every
$\varepsilon>0$ and every finite dimensional normed space $X$ there is a subspace $H$ of $X$ such that
$$
\operatorname{dim}(H)\geq c\varepsilon^{2}\ln(\operatorname{dim}(X))
\hskip 10mm\mbox{and}\hskip 10mm
d(H,\ell_{2}^{\operatorname{dim}(H)})\leq1+\varepsilon,
$$
where $\ell_{2}^{k}$ denotes the Euclidean space of dimension $k$.
Also this theorem is proved by the probabilistic method.
But even more is proved: not only there exists such a subspace, but actually most of the subspaces satisfy these specifications. This leads to an interesting paradoxon: whenever one tries
to select a subspace that is almost an Euclidean subspace one fails, but
if one chooses a subspace randomly, then with high probability one chooses a subspace that is almost Euclidean.
\par
As it turns out, Dvoretzky's theorem is related to a fundamental problem in Quantum Information Theory. The goal is to determine the capacity of a quantum channel to transmit classical information and a question that naturally arose in this context, and which has been one of the major open problems for more than a decade, asks about the additivity of the so-called  $\chi$-quantity for any pair of quantum channels. A more tractable, but equivalent question (see \cite{Sh2004}), concerns the additivity of the minimal output von Neumann entropy of quantum channels. Eventually, this reformulation allowed Hastings \cite{Hast} to construct his famous counterexample.
Shortly after, it was discovered by Aubrun, Szarek, and Werner in \cite{ASW1, ASW2} that the existence of this counterexample is a consequence of Dvoretzky's theorem.
\par
A classical result in convex geometry and the local theory of Banach spaces is John's theorem \cite{J1948} (see also \cite{B1992}) on maximal volume ellipsoids in convex bodies, which shows that, for any $n$-dimensional normed space $X$,
\[
d(X,\ell_2^n) \leq \sqrt{n}.
\]
The latter implies that $d(X,Y)\leq n$ for any pair $X$, $Y$ of $n$-dimensional normed spaces. It might come as a surprise that $d(\ell_1^n,\ell_\infty^n)$ is only of order $\sqrt{n}$ and this raises the question whether or not a pair of spaces can be constructed that has Banach-Mazur distance of order $n$. The following major achievement is due to Gluskin \cite{Glu1981}: there is a constant $c\in(0,\infty)$ such that for all $n\in\mathbb N$
there are $n$-dimensional normed spaces $X$ and $Y$ such that
$$
d(X,Y)\geq c\, n.
$$
Again, this result is proved by the probabilistic method and involves certain symmetric random polytopes acting as unit balls of the random normed spaces.
\par

Finally, let us mention some recent results from information-based complexity (IBC) concerning the quality of random information in approximation problems compared to optimal information. A typical question in IBC is to approximate the solution of a linear problem based on $n\in\N$ pieces of information about the unknown problem instance. Usually, it is assumed that some kind of oracle is available which grants us this information at our request and we may call this oracle $n$ times to get $n$ well-chosen pieces of information,
trying to obtain optimal information about the problem instance. However, often this model is too idealistic and there might be no oracle at our disposal and
the information just comes in randomly. In \cite{HKNPU2021}, the authors studied the circumradius of the intersection of an $m$-dimensional ellipsoid~$\mathcal E$
with semi-axes $\sigma_1\geq\dots\geq \sigma_m$ with random subspaces of
codimension $n$, where $n$ can be much smaller than $m$, and proved that under certain assumptions on
$\sigma$, 
this random radius 
$\mathcal{R}_n=\mathcal{R}_n(\sigma)$ 
is of the same order as the minimal such 
radius $\sigma_{n+1}$ with high probability. In other situations $\mathcal{R}_n$ 
is close to the maximum~$\sigma_1$. The random variable $\mathcal{R}_n$ naturally corresponds to the worst-case error of the best algorithm based on random information for 
$L_2$-approximation of functions from a compactly embedded 
Hilbert space $H$ with unit ball $\mathcal E$. 
In particular, $\sigma_k$ is the $k$th largest singular value of the embedding 
$H\hookrightarrow L_2$. In this formulation, one may also consider the 
case $m=\infty$ and it was shown that 
random information underlies an $\ell_2$-dichotomy in that it behaves differently depending 
on whether $\sigma \in \ell_2$ or not. We also refer the reader to subsequent works \cite{HPS2021, KU2021} and the survey \cite{HKNPU19-survey}.

In this paper we give an overview of some results concerning
best and random approximation of convex bodies by polytopes and how both are linked.
As we shall see random approximation is almost as good as
best approximation.
\par
The accuracy of approximation of a convex body by a polytope is interesting in itself, but it is also
relevant in many applications, for instance in
computervision (\cite{SaTa1}, \cite{SaTa2}), tomography \cite{Ga}, geometric
algorithms \cite{E}. 

\section{Notation \& Preliminaries}

We shall briefly set out the notation and some basic concepts used in this paper.  By $B_{2}^{n}(x,r)$ we denote the closed Euclidean ball with center $x\in\R^n$ and radius $r\in(0,\infty)$. The Euclidean norm on $\mathbb R^{n}$ is denoted by
$\|\ \|_{2}$ and we write $S^{n-1}:=\{x\in\R^n\,:\, \|x\|_2=1 \}$ for the Euclidean unit sphere. The standard inner product on $\R^n$ is denote by $\langle \cdot,\cdot \rangle$. A convex body $K$ in $\mathbb R^{n}$ is a compact convex set with non-empty interior. For such a body, the surface measure on its boundary, we write $\mu_{\partial K}$,
is the restriction of the $(n-1)$-dimensional Hausdorff measure
$\mathcal H^{n-1}$ to $\partial K$. For $x \in \partial K$ the normal at x to $\partial K$ is
denoted by $N_{K}(x)$. $N_{K}(x)$ is almost everywhere unique. We shall denote by
$\kappa_{K}(x)$  the Gau{\ss}-Kronecker
curvature of $\partial K$ at $x$.
\par
The Hausdorff distance between two convex bodies $C$ and $K$ in 
$\mathbb R^{n}$ is
$$
d_{H}(C,K):=\inf\{t>0\,:\,C\subseteq K+tB_{2}^{n}\hskip 1mm\mbox{and}\hskip 1mm   K\subseteq C+tB_{2}^{n}\}.
$$
The symmetric difference distance or Nikodym metric is defined as
$$
d_{S}(C,K):=\operatorname{vol}_{n}(C\triangle K)
=\operatorname{vol}_{n}(C\setminus K)+\operatorname{vol}_{n}( K\setminus C)
=\operatorname{vol}_{n}(C\cup K)-\operatorname{vol}_{n}(C\cap K),
$$
where $\operatorname{vol}_{n}(\cdot)$ refers to the $n$-dimensional Lebesgue measure.
We focus in this paper on the symmetric difference metric.
\par
For points $x_1, \dots x_N \in\Bbb R^{n}$ we denote by
$$
[x_1, \dots x_N]=\left\{\lambda_1 x_1 +\dots +\lambda _N x_N
\,:\, \forall1\leq i \leq N: 0\leq
\lambda_i\leq 1
\hskip 1mm\mbox{and}\hskip 1mm\sum_{i=1}^{N}\lambda_{i}=1\right\}
$$
the convex hull of these points. 
For two sequences $(a_n)_{n\in\N}, (b_n)_{n\in\N}\in\R^\N$, we write $a_n\sim b_n$ if $\lim_{n\to\infty} a_n/b_n=1$. Moreover, if the sequences are non-negative, then we use the notation $a_n\lesssim b_n$ to indicate that there exists a constant $C\in(0,\infty)$ such that, for all $n\in\N$, $a_n \leq C b_n$. We write $a_n \approx b_n$ if both $a_n \lesssim b_n$ and $b_n \lesssim a_n$.

\section{Best approximation of a convex body by polytopes contained in it}
Bronshte\u{\i}n and Ivanov \cite{BI} and Dudley \cite{Dud}
proved independently that for every convex body $C$ in $\mathbb R^{n}$
there is a constant $a_{C}:=a_{C}(n)\in(0,\infty)$ such that for all $N\in\mathbb N$ there exists 
a polytope $P_{N}$ with at most $N$ vertices such that
\begin{equation}\label{ApproxHaus1}
d_{H}(C,P_{N})\leq \frac{a_{C}}{N^{\frac{2}{n-1}}},
\end{equation}
i.e., with respect to the Hausdorff distance, convex bodies can be approximated arbitrarily well by polytopes (see also \cite{SW1981}). Note that the dependence in $N$ best possible.
Since 
$$
d_{S}(C,K)\leq d_{H}(C,K)\,\mathcal H_{n-1}(\partial (C\cup K)),
$$
the inequality carries over to the symmetric difference metric.
Nearly two decades later, the estimate in \eqref{ApproxHaus1} has been made more precise by Gordon, Meyer, and Reisner
\cite{GMR1}, showing the following: there is a constant $c\in(0,\infty)$ such that for all convex bodies
$K$ in $\mathbb R^{n}$ and all $N\in\mathbb N$ there is
a polytope $P_{N}$ having at most $N$ vertices such that
\begin{equation}\label{ApproxBall1-1}
d_{S}(K,P_{N})\leq\frac{c\cdot n}{N^{\frac{2}{n-1}}}
\operatorname{vol}_{n}(K),
\end{equation}
i.e., the constant $a_c$ has been specified. Of course, it remains to understand the constant $c$ and also whether the dependence in terms of volume of the body $K$ is optimal or if this dependence can be replaced by a quantity that is related to the surface structure of $K$. We will come back to this question later.

Macbeath \cite{Macb} observed that for every convex body $C$ in
$\mathbb R^{n}$ with 
$\operatorname{vol}_{n}(C)=\operatorname{vol}_{n}(B_{2}^{n})$,
for every $N\in\mathbb N$
and every polytope $P_{N}$ contained in $B_{2}^{n}$ having at most
$N$ vertices there is a polytope $Q_{N}$ contained in $C$ having at most $N$ vertices and
$$
\operatorname{vol}_{n}(P_{N})
\leq \operatorname{vol}_{n}(Q_{N}),
\hskip 10mm\mbox{i.e.}\hskip 10mm
d_{S}(Q_{N},C)\leq d_{S}(P_{N},B_{2}^{n}).
$$
Therefore, to prove the inequality (\ref{ApproxBall1-1})
it is enough to show it for the Euclidean ball (and polytopes contained in the Euclidean ball). By a result of Kabatjanskii and Levenstein \cite{KabLev}
for any given angle $\phi$ there are $\xi_{1},\dots,\xi_{N}\in S^{n-1}$
such that
\begin{eqnarray*}
\cos\phi\geq\langle \xi_{i},\xi_{j}\rangle
\hskip 20mm i\ne j
\\
\forall x\in S^{n-1}\exists i:
\hskip 10mm \cos\phi\leq\langle x,\xi_{i}\rangle
\\
N\leq (1-\cos\phi)^{-\frac{n-1}{2}}2^{0.901(n-1)}.
\end{eqnarray*}
Choosing
$P_{N}=[\xi_{1},\dots,\xi_{N}]$, we obtain
\begin{equation}\label{BestApproxAbove1}
 \operatorname{vol}_{n}(B_{2}^{n})- \operatorname{vol}_{n}(P_{N})
 \leq 2^{1.802-2}
 \frac{ n\cdot\operatorname{vol}_{n}(B_{2}^{n})}{N^{\frac{2}{n-1}}}.
\end{equation}
On the other hand, Gordon, Reisner and Sch\"utt \cite{GRS1,GRS2}
showed the following.

\begin{proposition}\label{BestApproxBelow1}
There are constants $a,b\in(0,\infty)$ such that for every $n\geq2$
and every $N\geq(b\cdot n)^{\frac{n-1}{2}}$, and every polytope $P_{N}$
contained in the Euclidean unit ball having at most $N$ vertices,
\begin{equation}\label{BestApproxBelow1-1}
\frac{a\cdot n\cdot\operatorname{vol}_{n}(B_{2}^{n})}
{N^{\frac{2}{n-1}}}
\leq\operatorname{vol}_{n}(B_{2}^{n})-\operatorname{vol}_{n}(P_{N}).
\end{equation}
\end{proposition}

Here the condition $N\geq(b\cdot n)^{\frac{n-1}{2}}$ is needed to ensure the specific
constant in (\ref{BestApproxBelow1-1}). If we drop the assumption $N\geq(b\cdot n)^{\frac{n-1}{2}}$, then the estimate becomes
$$
\frac{a\cdot\operatorname{vol}_{n}(B_{2}^{n})}
{N^{\frac{2}{n-1}}}
\leq\operatorname{vol}_{n}(B_{2}^{n})-\operatorname{vol}_{n}(P_{N}).
$$
We believe though that inequality (\ref{BestApproxBelow1-1}) holds without this assumption.

Let us see if more precise estimates for best approximation
of a convex body by polytopes can be obtained. The following result relates the symmetric difference between a convex body with sufficiently smooth boundary and a polytope contained in it with the so-called Delone triangulation constant and the affine surface area of $K$.

\begin{theorem}\label{GruberBest}
For all $n\in\mathbb N$ with $n\geq2$ and for all convex bodies 
$K$ in $\Bbb R^{n}$ with $C^{2}$-boundary 
and everywhere positive curvature, we have (asymptotically in $N$)
\begin{align*}
& \inf \{d_{S}(K,P_{N})\,:\,P_{N} \subseteq K\ \mbox{and}  
\hskip 2mm  P_{N} \mbox{ is a polytope with  at most N vertices} \} \cr
& \sim \tfrac{1}{2} \operatorname{del}_{n-1}
\left(\int_{\partial K} \kappa_{K}(x)^{\frac{1}{n+1}}
d\mu_{\partial K}(x)\right)^{\frac{n+1}{n-1}} 
\left(\frac{1}{N}\right)^{\frac{2}{n-1}},
\end{align*}
where $\operatorname{del}_{n-1}$ is a constant that is related to the Delone triangulation
and depends only on the dimension $n$. 
\end{theorem}

The integral expression 
\[
\int_{\partial K} \kappa_{K}(x)^{\frac{1}{n+1}}
d\mu_{\partial K}(x)
\]
is a fundamental quantity in convex geometry and called the affine surface area of $K$.

Theorem \ref{GruberBest} was proved by
McClure and Vitale \cite{McCVitale} in dimension 2  and by Gruber \cite{Gruber1,Gruber2} in higher dimensions.

\par
By (\ref{BestApproxAbove1}) and (\ref{BestApproxBelow1-1}) it follows that the constant del$_{n-1}$ is of the order of $n$, which
means that there are
numerical constants $a,b\in(0,\infty)$ such that, for all $n\in\Bbb N$,
\begin{equation}\label{delunif1}
a\cdot n\leq \mbox{del}_{n-1}\leq b\cdot n.
\end{equation}
In two papers by Mankiewicz and
Sch\"utt the constant del$_{n-1}$ has been better estimated \cite{MankS1, MankS2}.
It was shown there that
\begin{equation}\label{MSdel1}
\frac{\tfrac{n-1}{n+1}}
{\operatorname{vol}_{n-1}(B_{2}^{n-1})^{\frac{2}{n-1}}}
\leq \operatorname{del}_{n-1} \leq
(1+\tfrac{c\ln n}{n})
\frac{\tfrac{n-1}{n+1}}
{\operatorname{vol}_{n-1}(B_{2}^{n-1})^{\frac{2}{n-1}}},
\end{equation}
where $c\in(0,\infty)$ is a numerical constant. It follows from Stirling's formula that
$$
\lim_{n\to\infty}\frac{\mbox{del}_{n-1}}{n}
=\frac{1}{2\pi e}=0.0585498....
$$
The right-hand inequality of (\ref{MSdel1}) is proved using a result by M\"uller \cite{Mue1}
on random polytopes of the Euclidean ball. It is a special case of Theorem \ref{TheoremSW1}.
Altogether we have quite precise estimates for best approximation 
of convex bodies by polytopes if the numbers of vertices of the polytopes are big.
\par
We try now to give a more precise estimate for best approximation
for the whole range of the numbers $N$ of vertices. While the asymptotic
estimate for large $N$ involves  the affine surface area, we involve in general
the convex floating body. We explain why this is quite natural and, of course, what the convex floating body is.
\par 
Let $K$ be a convex body in $\mathbb R^{n}$ and $t\geq0$.
The convex floating body $K_{t}$ of $K$ is the intersection of
all halfspaces $H^{+}$ whose defining hyperplanes $H$ cut off a set of volume
$\delta$ from $K$, i.e., 
$$
K_{\delta}
=\bigcap_{\operatorname{vol}_{n}
(H^{-}\cap K)=\delta}H^{+}.
$$
Let us introduce the notion of generalized Gau\ss-Kronecker curvature. A convex function
$f:X \rightarrow \Bbb R, X \subseteq \Bbb R^d$ is called twice differentiable at $x_0$
in a generalized sense if there are a linear map $d^{2}f(x_0) \in L(\Bbb R^d)$ and a 
neighborhood $U(x_0)$ so that we have for all $x \in U(x_0)$ and all subdifferentials
df(x)
$$   
\parallel df(x)-df(x_0)-(d^{2}f(x_0))(x-x_0) \parallel_{2}
     \leq \Theta(\parallel x-x_0 \parallel_{2} ) \parallel x-x_0 \parallel_{2}      ,
$$
where $\lim_{t \to 0} \Theta(t)=\Theta(0)=0$ and where $\Theta$ is a montone
function. $d^{2}f(x_0)$ is symmetric and positive semidefinite. If f(0)=0
and df(0)=0 then the ellipsoid or elliptical cylinder

$$ x^{t}d^{2}f(0)x=1    $$

is called the indicatrix of Dupin at 0. The general case is reduced to the case
f(0)=0 and df(0)=0 by an affine transform. The eigenvalues of $d^{2}f(0)$ are
called the principal curvatures and their product the Gau\ss-Kronecker curvature
$\kappa(0)$. 
\par
Busemann and Feller \cite{BuFe} proved for convex functions mapping from $\mathbb R^{2}$ to $\mathbb R$ that they have almost everywhere a generalized Gauss-Kronecker curvature.
Aleksandrov \cite{Aleks} extended their result to higher dimensions.

\begin{theorem}\label{SWFloat1}\cite{SchuWer}
For all $n\in\mathbb N$ with $n\geq2$ and all convex bodies $K$
in $\mathbb R^{n}$,
\begin{equation}\label{FloatAffS1}
\lim_{\delta \to 0} \frac{\operatorname{vol}_{n}(K)-\operatorname{vol}_{n}(K_{\delta})}
{\delta^{\frac{2}{n+1}}}
=
\frac{1}{2}\left(\frac{n+1}{\operatorname{vol}_{n-1}(B_{2}^{n-1})}\right)^{\frac{2}{n+1}}
\int_{\partial K}\kappa_{K}(x)^{\frac{1}{n+1}}d\mu_{\partial K}(x).
\end{equation}
\end{theorem}

The expression 
$$
\int_{\partial K}\kappa_{K}(x)^{\frac{1}{n+1}}d\mu_{\partial K}(x)
$$
is called affine surface area of $K$.
It is generally conjectured that Theorem \ref{GruberBest} holds for arbitrary convex bodies
if we substitute the classical Gau{\ss}-Kronecker curvature by the generalized one.

Since the affine surface area is an integral part of the formula for best approximation in the case of a large number of vertices and since the affine surface area is given by (\ref{FloatAffS1}) one can speculate whether
in the case of general numbers $N$ of vertices the floating body
is involved. In fact, this is the case.

Let us describe the floating body algorithm.
We are choosing the vertices $x_{1},\dots,x_{N}\in\partial K$
of the polytope $P_{N}$. $x_{1}$ is chosen arbitrarily.
Having chosen $x_{1},\dots,x_{k}$, we choose a support hyperplane $H_{k+1}$
to $K_{\delta}$
such that 
$$
K_{\delta}\subseteq H_{k+1}^{+}
\hskip 10mm\mbox{and}\hskip 10mm
x_{1},\dots,x_{k}\in H_{k+1}^{+}.
$$
If there is no such hyperplane the algorithm stops.
We choose $x_{k+1}\in\partial K\cap H_{k+1}^{-}$ such that the distance
of $x_{k+1}$ to $H_{k+1}$ is maximal. This means that
$x_{k+1}$ is contained in the support hyperplane $H$ to $K$
that is parallel to $H_{k+1}$ and such that $H\subseteq H_{k+1}^{-}$.
\vskip 5mm

\begin{theorem}\cite{Schu3}\label{PolyApproxFloat1}
Let $K$ be a convex body in $\mathbb R^{n}$. Then, for all $\delta$
with $0< \delta\leq\frac{1}{4e^{4}}\operatorname{vol}_{n}(K)$,
there exist $N\in\mathbb N$ with 
\begin{equation}\label{PolyApproxFloat1-1}
N
\leq
e^{16n}\frac{\operatorname{vol}_{n}(K\setminus K_{\delta})}
{\operatorname{vol}_{n}(B_{2}^{n})\delta}
\end{equation}
and a polytope $P_{N}$ that has at most $N$ vertices 
such that
$$
K_{\delta}\subseteq P_{N}\subseteq K.
$$
\end{theorem}

We want to see how good this algorithm is. Therefore,
we determine what Theorem \ref{PolyApproxFloat1}
gives in the case of a smooth body and in the case of a polytope.
First the case of a smooth body.

\begin{corollary} \label{PolyApproxFloat2}
Let $K$ be a convex body in $\mathbb R^{n}$. For every sufficiently large
$N\in\mathbb N$ there exists $\delta_{N}\in(0,\infty)$ with
$$
N\leq e^{16n}\frac{\operatorname{vol}_{n}(K\setminus K_{\delta_{N}})}
{\delta_{N}\operatorname{vol}_{n}(B_{2}^{n})}
$$
and a polytope $P_{N}$ with at most $N$ vertices that is constructed by the floating body algorithm
for the floating body $K_{\delta_{N}}$ such that
\begin{equation}\label{PolyApproxFloat2-1}
\limsup_{N\to\infty}\frac{d_{s}(K,P_{N})}{N^{-\frac{2}{n-1}}}
\leq c n^{2}\left(\int_{\partial K}\kappa^{\frac{1}{n+1}}d\mu_{\partial K}
\right)^{\frac{n+1}{n-1}}.
\end{equation}
\end{corollary}

The main difference between the optimal result and (\ref{PolyApproxFloat2-1}) is that $n^{2}$ appears as a factor
on the right-hand side of (\ref{PolyApproxFloat2-1}) and not $n$.
It still an open question whether or not an improvement of the argument gives the same estimate but with $n$ instead of $n^{2}$.
\vskip 3mm

\noindent
{\bf Proof of Corollary \ref{PolyApproxFloat2}.} By (\ref{FloatAffS1})
$$
\lim_{\delta \to \infty} \frac{\operatorname{vol}_{n}(K)-\operatorname{vol}_{n}(K_{\delta})}{\delta^{\frac{2}{n+1}}}=
\frac{1}{2}\left(\frac{n+1}{\operatorname{vol}_{n-1}(B_{2}^{n-1})}\right)^{\frac{2}{n+1}}
\int_{\partial K}\kappa(x)^{\frac{1}{n+1}}d\mu_{\partial K}(x).
$$
Therefore,
$$
\operatorname{vol}(K\setminus K_{\delta})
\sim \delta^{\frac{2}{n+1}}\frac{1}{2}
\left(\frac{n+1}{\operatorname{vol}_{n-1}(B_{2}^{n-1})}\right)^{\frac{2}{n+1}}
\int_{\partial K}\kappa^{\frac{1}{n+1}}d\mu_{\partial K}
$$
and 
$$
e^{16n}\frac{\operatorname{vol}_{n}(K\setminus K_{\delta})}
{\delta\operatorname{vol}_{n}(B_{2}^{n})}
\approx \delta^{-\frac{n-1}{n+1}}e^{16n}
\frac{1}{2}
\frac{(n+1)^{\frac{2}{n+1}}}{\left(\operatorname{vol}_{n-1}(B_{2}^{n-1})
\right)^{\frac{n+3}{n+1}}}
\int_{\partial K}\kappa^{\frac{1}{n+1}}d\mu_{\partial K}
=\delta^{-\frac{n-1}{n+1}}c(n,K).
$$
Therefore, there is $N_{0}\in\N$ such that for all $N\geq N_{0}$
there is $\delta_{N}\in(0,\infty)$ with
$$
N\leq e^{16n}\frac{\operatorname{vol}_{n}(K\setminus K_{\delta_{N}})}
{\delta_{N}\operatorname{vol}_{n}(B_{2}^{n})}
<N+1.
$$
It follows that 
\begin{eqnarray*}
&&\frac{\operatorname{vol}_{n}(K\setminus P_{N})}
{\frac{1}{N^{\frac{2}{n-1}}}}
\leq\frac{\operatorname{vol}_{n}(K\setminus K_{\delta_{N}})}
{\frac{1}{N^{\frac{2}{n-1}}}}  
\leq\frac{\operatorname{vol}_{n}(K\setminus K_{\delta_{N}})}
{( \frac{\delta_{N}\operatorname{vol}_{n}(B_{2}^{n})}{e^{16n}\operatorname{vol}_{n}(K\setminus K_{\delta_{N}})}
)^{\frac{2}{n-1}}} \\
&&=\frac{e^{\frac{32n}{n-1}}}{(\operatorname{vol}_{n}(B_{2}^{n}))^{\frac{2}{n-1}}}
\frac{(\operatorname{vol}_{n}(K\setminus K_{\delta_{N}}))^{\frac{n+1}{n-1}}}
{ \delta_{N}^{\frac{2}{n-1}}}  
\approx e^{32}n\left(
\frac{\operatorname{vol}_{n}(K\setminus K_{\delta_{N}})}{ \delta_{N}^{\frac{2}{n+1}}} \right)^{\frac{n+1}{n-1}}  \\
&&\approx e^{32}n\left(
\frac{1}{2}\left(\frac{n+1}{\operatorname{vol}_{n-1}(B_{2}^{n-1})}\right)^{\frac{2}{n+1}}
\int_{\partial K}\kappa(x)^{\frac{1}{n+1}}d\mu_{\partial K}(x) \right)^{\frac{n+1}{n-1}}   \\
&&\approx e^{32}n^{2}\left(
\int_{\partial K}\kappa(x)^{\frac{1}{n+1}}d\mu_{\partial K}(x) \right)^{\frac{n+1}{n-1}},
\end{eqnarray*}
which completes the proof.
$\Box$
\vskip 3mm

On the other end of the spectrum of convex bodies are
the polytopes.  Of course, best approximation of a polytope
by polytopes is a trivial task once the number of vertices one may choose 
is equal to or greater than the number of vertices of the polytope.
We cannot expect that Theorem \ref{PolyApproxFloat1} gives us this. But,
what does Theorem \ref{PolyApproxFloat1} give in the case of a polytope?

Let us denote the set of all polytopes in $\mathbb R^{n}$ by
$\mathcal P^{n}$.
The extreme points or vertices of a polytope $P$ shall be denoted by $\operatorname{ext}(P)$. 
\par
A $n+1$-tuple $(f_{0},\dots,f_{n})$ such that $f_{i}$ is an $i$-dimensional face of
a polytope $P$ for all $i=0,1,\dots,n$ and such that
$$
f_{0}\subset f_{1}\subset\cdots\subset f_{n}
$$
is called a flag or tower of $P$. We denote the set of all flags of $P$ by
$\operatorname{\mathcal Fl}(P)$. The number of all flags of a polytope $P$
is denoted by $\operatorname{flag}(P)$. Clearly, $f_{n}=P$ and $f_{0}$ is a vertex of $P$.
\par
We establish two recurrence formulae for $\operatorname{flag}(P)$.
We define
$\phi_{n}:\mathcal P^{n}\to\mathbb R$ by
$$
\phi_{1}(P)=2
$$
and
$$
\phi_{n}(P)
=\sum_{x\in\operatorname{ext}(P)}\phi_{n-1}(P\cap H_{x}),
$$
where $H_{x}$ is a hyperplane that separates $x\in\R^n$ strictly
from all other extreme points of $P$. 
\par
Moreover, we define $\psi_{n}:\mathcal P^{n}\to\mathbb R$ by
$$
\psi_{1}(P)=2
$$
and for $n\geq2$ by
$$
\psi_{n}(P)
=\sum_{F\in\operatorname{fac}_{n-1}(P)}
\psi_{n-1}(F),
$$
where $\operatorname{fac}_{n-1}(P)$ is the set of all
$n-1$-dimensional faces of $P$.
\par
It can be shown that
$$
\phi_{n}(P)=\psi_{n}(P)=\operatorname{flag}(P).
$$
We can further prove that for all polytopes with $0$ as an interior point
and $P^{\circ}$ its dual/polar polytope 
$$
\operatorname{flag}(P)=\operatorname{flag}(P^{\circ}),
$$
where $P^{\circ}:= \{y\in\R^n\,:\,\forall x\in P:\,\, \langle x,y \rangle \leq 1 \}$. 
It is not difficult to show for the cube $C^{n}$ and the simplex $\Delta^{n}$
in $\mathbb R^{n}$
$$
\operatorname{flag}(\Delta^{n})=n!
\hskip 20mm
\operatorname{flag}(C^{n})=\operatorname{flag}((C^{n})^{\circ})=2^{n}n!.
$$
For the floating bodies of polytopes we have

\begin{theorem}\cite{Schu1}\label{PolyRandomPoly}
\label{PolyFloatBod1}
Let $P$ be a convex polytope with nonempty interior
in $\mathbb R^{n}$. Then
\begin{equation}\label{PolyFloatBod1-1}
\lim_{\delta\to0}\frac{\operatorname{vol}_{n}(P)
-\operatorname{vol}_{n}(P_{\delta})}{\delta\left(\ln\frac{1}{\delta}\right)^{n-1}}
=\frac{\operatorname{flag}_{n}(P)}{n!n^{n-1}}.
\end{equation}
\end{theorem}

Theorem \ref{PolyFloatBod1} was proved by Sch\"utt \cite{Schu1}. Note that B\'ar\'any and Larman \cite{BL} showed that for all polytopes with an interior point in $\mathbb R^{n}$, we have
$$
\operatorname{vol}_{n}(P\setminus P_{\delta})
\leq c \delta\left(\ln\frac{1}{\delta}\right)^{n-1}.
$$

The previous theorem implies fast convergence of the approximation of a given polytope by a polytope with at most $N$ vertices. Indeed, by (\ref{PolyFloatBod1-1}), for sufficiently small $\delta\in(0,\infty)$, we have
$$
\operatorname{vol}_{n}(P)-\operatorname{vol}_{n}(P_{\delta})
\sim\frac{\operatorname{flag}_{n}(P)}{n!n^{n-1}}
\delta\left(\ln\frac{1}{\delta}\right)^{n-1}.
$$
Therefore, for sufficiently small $\delta\in(0,\infty)$,
\begin{eqnarray*}
\operatorname{vol}_{n}(P)-\operatorname{vol}_{n}(P_{N})
\leq\operatorname{vol}_{n}(P)-\operatorname{vol}_{n}(P_{\delta})
\sim\frac{\operatorname{flag}_{n}(P)}{n!n^{n-1}}
\delta\left(\ln\frac{1}{\delta}\right)^{n-1}.
\end{eqnarray*}
By (\ref{PolyApproxFloat1-1})
$$
N\leq e^{16n}\frac{\operatorname{vol}_{n}(P)-\operatorname{vol}_{n}(P_{\delta})}
{\operatorname{vol}_{n}(B_{2}^{n})\delta}
\sim e^{16n}\frac{\operatorname{flag}_{n}(P)}
{n!n^{n-1}\operatorname{vol}_{n}(B_{2}^{n})}\left(\ln\frac{1}{\delta}\right)^{n-1}.
$$
This implies
$$
N^{\frac{1}{n-1}}e^{-32}\left(\frac{n!n^{n-1}\operatorname{vol}_{n}(B_{2}^{n})}{\operatorname{flag}_{n}(P)}\right)^{\frac{1}{n-1}}
\leq\ln\frac{1}{\delta}
$$
and
$$
\delta\leq e^{-N^{\frac{1}{n-1}}}\exp\left(-e^{-32}\left(\frac{n!n^{n-1}\operatorname{vol}_{n}(B_{2}^{n})}{\operatorname{flag}_{n}(P)}\right)^{\frac{1}{n-1}}\right)
\sim e^{-N^{\frac{1}{n-1}}}\exp\left(-\frac{e^{-32}n^{\frac{3}{2}}}{(\operatorname{flag}_{n}(P))^{\frac{1}{n-1}}}\right).
$$
Since $\delta\left(\ln\frac{1}{\delta}\right)^{n-1}$ is an increasing function
for $\delta$ with $0\leq\delta\leq e^{-n+1}$, we obtain
\begin{eqnarray*}
&&\operatorname{vol}_{n}(P)-\operatorname{vol}_{n}(P_{N})
\\
&&\lesssim\frac{\operatorname{flag}_{n}(P)}{n!n^{n-1}}
e^{-N^{\frac{1}{n-1}}}\exp\left(-\frac{e^{-32}n^{\frac{3}{2}}}{(\operatorname{flag}_{n}(P))^{\frac{1}{n-1}}}\right)
\left(N^{\frac{1}{n-1}}+\frac{e^{-32}n^{\frac{3}{2}}}{(\operatorname{flag}_{n}(P))^{\frac{1}{n-1}}}\right)^{n-1}.
\end{eqnarray*}
This means that the right-hand side expression is of the order
$$
\frac{N}{e^{N^{\frac{1}{n-1}}}},
$$
which is clearly decreasing rapidly to $0$ for $N\to\infty$.

\section{Random approximation of convex bodies by polytopes}

We now turn to the topic of random approximation.
Let $K$ be a convex body in $\Bbb R^n$. A random polytope in $K$ is the convex
hull of finitely many points in $K$ that are chosen at random with respect
to a probability measure on $K$. First we consider the natural choice of the normalized Lebesgue measure on $K$, i.e., we consider the uniform distribution on $K$.
For a fixed number $N$ of points we are interested in the expected volume
of that part of $K$ that is not contained in the convex hull $[x_1,.....,
x_N]$ of the chosen points. Let us denote
\begin{equation}\label{ExpVol1}
 \Bbb E(K,N)= \int_{K \times \cdots \times K} \operatorname{vol}_n([x_1,...,x_N])\,d\Bbb P(x_1,...x_N)   ,
 \end{equation}
where $\Bbb P$ is the $N$-fold product of the normalized Lebesgue measure on $K$.
We are interested in the asymptotic behavior of
\begin{equation}\label{ExpVol2}
\operatorname{vol}_n(K)-\Bbb E(K,N)= \int_{K \times \cdots \times K} \operatorname{vol}_n(K \setminus
[x_1,....,x_N]) \,d\Bbb P(x_1,...,x_N)    . 
\end{equation}

\vskip 1cm

\begin{theorem}\label{ExpRandPoly}
Let K be a convex body in $\Bbb R^n$. Then we have
\begin{equation}\label{ExpRandPoly-1}
c(n)\lim_{N \to \infty} \frac {\operatorname{vol}_n(K)-\Bbb E(K,N)}{\left(\frac{\operatorname{vol}_n(K)}{N}\right)^{\frac{2}{n+1}}}
=\int_{\partial K} \kappa_{K}^{\frac{1}{n+1}}d\mu_{\partial K} ,           
\end{equation}
where $\kappa_{K} (x)$ is the generalized Gau\ss-Kronecker curvature and
$$ 
c(n)=2\left(\frac{\operatorname{vol}_{n-1}(B_2^{n-1})}{n+1}\right)^{\frac {2}{n+1}} \frac{(n+3)(n+1)!}
    {(n^2+n+2)(n^2+1)\Gamma(\frac{n^2+1}{n+1})}            .                        
    $$
\end{theorem}

R\'{e}nyi and Sulanke  \cite{RS1, RS2} proved (\ref{ExpRandPoly-1})
for smooth convex bodies in $\mathbb R^{2}$.
Wieacker \cite{Wie} extended their results to higher dimensions
for the Euclidean ball.
Schneider and Wieacker \cite{SchWi1}
extended the results to higher dimensions for the mean width
instead of the difference volume.
 It has been solved
by B\'ar\'any [B] for convex bodies with $C^3$ boundary and everywhere positive curvature. The result
 for arbitrary convex bodies had been shown in \cite{Schu2} (see also \cite{BHH}).
\par

By choosing the points randomly from the convex body $K$, we get a polytope $P_{N}$
with at most $N$ vertices and symmetric difference distance to $K$ less than
$$
\left(\frac{\operatorname{vol}_n(K)}{N}\right)^{\frac{2}{n+1}}
\frac{1}{2}\left(\frac{n+1}{\operatorname{vol}_{n-1}(B_2^{n-1})}\right)^{\frac {2}{n+1}} \frac
    {(n^2+n+2)(n^2+1)\Gamma(\frac{n^2+1}{n+1})} {(n+3)(n+1)!}
\int_{\partial K} \kappa_{K}^{\frac{1}{n+1}}d\mu_{\partial K}   ,
$$
while the optimal order is
$$
\tfrac{1}{2} \mbox{del}_{n-1}
\left(\int_{\partial K} \kappa_{K}(x)^{\frac{1}{n+1}}
d\mu_{\partial K}(x)\right)^{\frac{n+1}{n-1}} 
\left(\frac{1}{N}\right)^{\frac{2}{n-1}}.
$$
The random approximation is proportinal to 
$\left(\frac{1}{N}\right)^{\frac{2}{n+1}}$, while the best approximation
is proportional to $\left(\frac{1}{N}\right)^{\frac{2}{n-1}}$ by (\ref{ApproxBall1-1}) and (\ref{BestApproxBelow1-1}).
On the other hand, as shown in \cite{WeilW1}, the expected number of vertices 
of a random polytope equals
$$
 c_{n}\int_{\partial K} \kappa_{K}(x)^{\frac{1}{n+1}}
d\mu_{\partial K}(x)\left(\frac{N}{\operatorname{vol}_{n}(K)}\right)^{\frac{n-1}{n+1}},
$$
so that the number of vertices of a random polytope is of the order
$N^{\frac{n-1}{n+1}}$. This suggests that there exists a random polytope with 
$N$ vertices whose symmetric difference distance is of the order
$(\frac{1}{N})^{\frac{2}{n-1}}$.
\par
It is naturally to expect a better order of random approximation if the points are chosen from the boundary $\partial K$.

\begin{theorem}\label{TheoremSW1}
Let K be a convex body in $\mathbb{R}^n$ such that there are $r,R\in(0,\infty)$ with
$0<r\leq R < \infty$ so that, for all $x\in\partial K$,
\begin{equation}\label{TheoremSW1-1}
B_{2}^{n}(x-rN_{\partial K}(x),r)\subseteq K
\subseteq B_{2}^{n}(x-RN_{\partial K}(x),R)
\end{equation}
and let $f:\partial K\rightarrow\Bbb R_{+}$ be a continuous, positive function
with
$
\int_{\partial K}f(x)d\mu_{\partial K}(x)=1.
$
Let $\Bbb P_{f}$ be the probability measure on $\partial K$ given by
$
d\Bbb P_{f}(x)=f(x)d\mu_{\partial K}(x).
$
Then we have
\begin{equation}\label{TheoremSW1-2}
\lim_{N \to \infty} \frac{\mbox{\rm vol}_n(K)-\Bbb{E}(f
,N)}{\left(\frac{1}{N}\right)^\frac{2}{n-1}}=c_{n}\int_{\partial K}
\frac{\kappa(x)^{\frac{1}{n-1}}}{f(x)^{\frac{2}{n-1}}}d\mu_{\partial K}(x),
\end{equation}
where $\kappa$ is the (generalized) Gau\ss-Kronecker curvature and
\begin{equation}\label{TheoremSW1-3}
c_n=
\frac{(n-1)^{\frac{n+1}{n-1}}\Gamma \left(n+1+\tfrac{2}{n-1}\right)}
{2(n+1)!(\operatorname{ vol}_{n-2}(\partial B_{2}^{n-1}))^{\frac{2}{n-1}}}.
\end{equation} 
The minimum on the right-hand side is attained for the normalized affine
surface area measure with density with respect to the surface area measure $\mu_{\partial K}$ 
\begin{equation}\label{TheoremSW1-4}
f_{\operatorname{as}}(x)=\frac{\kappa(x)^{\frac{1}{n+1}}}{\int_{\partial
K}\kappa(x)^{\frac{1}{n+1}}d\mu_{\partial K}(x)}.
\end{equation}
\end{theorem}
\vskip 5mm

Theorem \ref{TheoremSW1} has been obtained by Sch\"utt and Werner in \cite{SW5} proved. M\"uller \cite{Mue1} proved this result for the Euclidean ball.
Reitzner \cite{Re4} proved this result for convex bodies with $C^{2}$-boundary and everywhere positive curvature.
\par

The condition (\ref{TheoremSW1-1})
is satisfied if $K$ has a $C^{2}$-boundary with everywhere positive curvature.
This follows from Blaschke's rolling theorem (\cite{Bl2}, p.118) and a generalization 
of it \cite{SchuWer}. Indeed, we can choose  
$$
r=\min_{x\in\partial K}\min_{1\leq i\leq n-1}r_{i}(x)
\qquad\text{and}\qquad
R=\max_{x\in\partial K}\max_{1\leq i\leq n-1}r_{i}(x),
$$
where $r_{i}(x)$ denotes the $i$-th principal curvature radius.
\par
We show now that the expression (\ref{TheoremSW1-2}) for any other measure given by a density $f$ with respect to $\mu_{\partial K}$
is greater than or equal to the one for (\ref{TheoremSW1-4}). Since 
$\int_{\partial K}f(x)d\mu_{\partial K}(x)=1$, we have
\begin{eqnarray*}
& & \left(\frac{1}{\mbox{vol}_{n-1}(\partial
K)}\int_{\partial K}\left|\frac{\kappa(x)}{f(x)^{2}}\right|^{\frac{1}{n-1}}
d\mu_{\partial K}(x)\right)^{\frac{1}{n+1}}   \\
& &  =\left(\frac{1}{\mbox{vol}_{n-1}(\partial
K)}\int_{\partial
K}\left|\left(\frac{\kappa(x)}{f(x)^{2}}\right)^{\frac{1}{n^{2}-1}}
\right|^{n+1}
d\mu_{\partial K}(x)\right)^{\frac{1}{n+1}}  
\int_{\partial K}f(x)d\mu_{\partial K}(x) \\
& &  =\left(\frac{1}{\mbox{vol}_{n-1}(\partial
K)}\int_{\partial
K}\left|\left(\frac{\kappa(x)}{f(x)^{2}}\right)^{\frac{1}{n^{2}-1}}
\right|^{n+1}
d\mu_{\partial K}(x)\right)^{\frac{1}{n+1}} \times  \\
& &
\hskip 5mm
\left(\mbox{vol}_{n-1}(\partial
K)\right)^{\frac{2}{n^{2}-1}} 
\left(\frac{1}{\mbox{vol}_{n-1}(\partial
K)}\int_{\partial
K}\left|f(x)^{\frac{2}{n^{2}-1}}\right|^{\frac{n^{2}-1}{2}}d\mu_{\partial
K}(x)\right)^{\frac{2}{n^{2}-1}}.
\end{eqnarray*}
We have $ \frac{1}{n+1}+\frac{2}{n^{2}-1}=\frac{1}{n-1}$ and 
we apply  H\"older's inequality to
get
\begin{eqnarray*}
& & \left(\frac{1}{\mbox{vol}_{n-1}(\partial
K)}\int_{\partial K}\left|\frac{\kappa(x)}{f(x)^{2}}\right|^{\frac{1}{n-1}}
d\mu_{\partial K}(x)\right)^{\frac{1}{n+1}}   \\
& &  \geq
\left(\frac{1}{\mbox{vol}_{n-1}(\partial
K)}\int_{\partial K}\kappa(x)^{\frac{1}{n+1}}d\mu_{\partial K}(x)
\right)^{\frac{1}{n-1}}\left(\mbox{vol}_{n-1}(\partial
K)\right)^{\frac{2}{n^{2}-1}},
\end{eqnarray*}
which gives us
$$
\int_{\partial K}\left|\frac{\kappa(x)}{f(x)^{2}}\right|^{\frac{1}{n-1}}
d\mu_{\partial K}(x)
\geq\left(\int_{\partial K}\kappa(x)^{\frac{1}{n+1}}d\mu_{\partial K}(x)
\right)^{\frac{n+1}{n-1}}.
$$
\par

Choosing points according to the affine surfaca area 
measure gives random polytopes of greatest possible volume.
Again, this is intuitively clear:  an optimal measure should put
more weight on points with higher curvature. Moreover, and this is a
crucial observation, if the optimal measure is unique, then it must be
affine invariant. This measure is affine invariant, i.e., for an affine, volume preserving map
$T$ and all measurable subsets $A$ of $\partial K$,
$$
\int_{A}\kappa_{\partial K}^{\frac{1}{n+1}}(x)d\mu_{\partial K}(x)
=\int_{T(A)}\kappa_{\partial T(K)}^{\frac{1}{n+1}}(x)d\mu_{\partial
T(K)}(x).
$$
There are not too many such measures and the affine
surface measure is the first that comes to ones mind. This measure
satisfies two necessary properties: it is affine invariant and it puts more
weight on points with greater curvature.
\par
The order of magnitude  for random
approximation for large $N$ is 
$$
\frac{(n-1)^{\frac{n+1}{n-1}}\Gamma \left(n+1+\tfrac{2}{n-1}\right)}
{2(n+1)!(\mbox{vol}_{n-2}(\partial B_{2}^{n-1}))^{\frac{2}{n-1}}}
\left(\int_{\partial K}
\kappa(x)^{\frac{1}{n+1}}d\mu_{\partial K}(x)\right)^{\frac{n+1}{n-1}}
\left(\frac{1}{N}\right)^\frac{2}{n-1}   , 
$$
while best approximation for large $N$ is of the order
$$
\tfrac{1}{2} \operatorname{del}_{n-1}
\left(\int_{\partial K} \kappa_{K}(x)^{\frac{1}{n+1}}
d\mu_{\partial K}(x)\right)^{\frac{n+1}{n-1}} 
\left(\frac{1}{N}\right)^{\frac{2}{n-1}}.
$$
The quotient between best and random approximation is
$$
\frac{(n-1)^{\frac{n+1}{n-1}}\Gamma \left(n+1+\tfrac{2}{n-1}\right)}
{2(n+1)!(\mbox{vol}_{n-2}(\partial B_{2}^{n-1}))^{\frac{2}{n-1}}}
\cdot\frac{2}{\mbox{del}_{n-1}}.
$$
Since $\tfrac{n-1}{n+1}\operatorname{vol}_{n-1}(B_{2}^{n-1})^{-\frac{2}{n-1}}
\leq \mbox{del}_{n-1} $, the quotient is less than
$$
\frac{(n-1)^{\frac{2}{n-1}}\Gamma \left(n+1+\tfrac{2}{n-1}\right)
\operatorname{vol}_{n-1}(B_{2}^{n-1})^{\frac{2}{n-1}}}
{n!(\operatorname{vol}_{n-2}(\partial B_{2}^{n-1}))^{\frac{2}{n-1}}}
=\frac{\Gamma\left(n+1+\frac{2}{n-1}\right)}{n!}.
$$
Since
$$
\lim_{x\to\infty}\frac{\Gamma(x+\alpha)}{x^{\alpha}\Gamma(x)}=1
$$
the quotient is asymptotically equal to
$$
(n+1)^{\frac{2}{n-1}}.
$$
Since the exponential function is convex, we have, for all $t\in[0,1]$,
$$
e^{t}\leq 1+(e-1)t.
$$
This implies
$$
(n+1)^{\frac{2}{n-1}}=\exp\left(\ln\left((n+1)^{\frac{2}{n-1}}\right)\right)
\leq1+2(e-1)\frac{\ln(n+1)}{n-1}
$$
and therefore the order of the quotient between random and best approximation
is less than 
$$
1+2(e-1)\frac{\ln(n+1)}{n-1}.
$$
We want to discuss two other measures that are of interest. 
The measure with the uniform density 
$$
f(x)=\frac{1}{\operatorname{vol}_{n-1}(\partial K)}
$$
with respect to $\mu_{\partial K}$
gives
$$
\lim_{N \to \infty} \frac{\operatorname{vol}_n(K)-\Bbb{E}(f
,N)}{\left(\frac{\operatorname{vol}_{n-1}(\partial
K)}{N}\right)^\frac{2}{n-1}}=c_{n}\int_{\partial K}
\kappa(x)^{\frac{1}{n-1}}d\mu_{\partial K}(x).
$$
 Let $K$ be a convex body and $\operatorname{cen}(K)$ its centroid. For any
 Borel set $A$ with
$A\subseteq \partial K$ the cone probability measure is defined by
$$
\mathbb P(A)=\frac{\operatorname{vol}_{n}([\operatorname{cen}(K),A])}{\operatorname{vol}_{n}(K)}.
$$
If the centroid is the origin, then the density is given by
$$
f(x)=\frac{\langle x,N_{\partial K}(x)\rangle}{\int_{\partial K}
\langle x,N_{\partial K}(x)\rangle d\mu_{\partial K}(x)}
$$
and the measure is invariant under linear, volume preserving maps.
We have 
$$
\frac{1}{n}\int_{\partial K}<x,N(x)>d\mu_{\partial K}(x)
=\operatorname{vol}_{n}(K)
$$
and thus
$$
f(x)=\frac{<x,N_{\partial K}(x)>}{n\ \mbox{vol}_{n}(K)}.
$$
We get
$$
\lim_{N \to \infty} \frac{\mbox{vol}_n(K)-\Bbb{E}(f
,N)}{\left(\frac{n\ \mbox{vol}_{n}(
K)}{N}\right)^\frac{2}{n-1}}=c_{n}\int_{\partial K}
\frac{\kappa(x)^{\frac{1}{n-1}}}{<x,N_{\partial K}(x)>^{\frac{2}{n-1}}}
d\mu_{\partial K}(x).
$$

\section{Approximation of convex bodies by polytopes without
cotainment condition}

So far we have considered the approximation of a convex body $K$
by polytopes that are contained in $K$. Now we want to drop the 
containment condition and find out whether or not the approximation improves. Intuitively, one expects an improvement in dependence on the dimension.
\par
Ludwig showed in \cite{Lud1} the following result.

\begin{theorem}\label{}\cite{Lud1}
Let $K$ be a convex body in $\mathbb R^{n}$ such that its
boundary is twice continuously differentiable and whose
curvature is everywhere strictly positive. Then
\begin{eqnarray*}
&&\lim_{N\to\infty}
\frac{\inf\{\operatorname{vol}_{n}(K\triangle P_{N})\,:\,P_{N} \hskip 1mm
\mbox{is a polytope with at most}
\hskip 1mm N\hskip 1mm\mbox{vertices}\}}
{N^{-\frac{2}{n-1}}}   \\
&&\hskip 30mm
=\frac{1}{2}\operatorname{ldel}_{n-1}\left(
\int_{\partial K}\kappa(x)^{\frac{1}{n+1}}d\mu_{\partial K}(x)\right)^{\frac{n+1}{n-1}}.
\end{eqnarray*}
\end{theorem}

The constant $\operatorname{ldel}_{n-1}$ is positive and depends only on
$n$.
Clearly, $\operatorname{ldel}_{n-1}\leq\operatorname{del}_{n-1}$
and by (\ref{delunif1}) it follows 
$\operatorname{ldel}_{n-1}\leq c\cdot n$. On the other hand, 
it has been shown in \cite{Boe} that for a polytope $P_{N}$
with at most $N$ vertices
$$
\operatorname{vol}_{n}(B_{2}^{n}\triangle P_{N})
\geq\frac{\operatorname{vol}_{n}(B_{2}^{n})}{67e^{2}\pi n}
\cdot\frac{1}{N^{\frac{2}{n-1}}}.
$$
Thus, between the upper and lower estimate for ldel$_{n-1}$
there is a gap of order $n^{2}$. In \cite{LSW}  Ludwig, Sch\"utt, and Werner 
 narrowed this gap by a factor $n$. It had been shown that ldel$_{n-1}\leq c$, where $c\in(0,\infty)$ is a numerical
constant.
\vskip 5mm

\begin{theorem}\label{vertex}\cite{LSW}
There is a constant $c\in(0,\infty)$ such that for every $n\in\mathbb N$
there is $N_{n}$ so that for
every
$N\geq N_{n}$ there is a polytope $P_{}$ in $\mathbb R^{n}$
with
$N$ vertices such that
\begin{equation}\label{improved bound}
\mbox{\rm vol}_{n}(B_{2}^{n}\triangle P_{})
\leq c\ \mbox{\rm vol}_{n}(B_{2}^{n})N^{-\frac{2}{n-1}}.
\end{equation}
\end{theorem}
Comparing the bound in \eqref{improved bound} with \eqref{ApproxBall1-1}, we see that the bound has been improved by a factor of the dimension. 
Therefore, there are constants $c_{1},c_{2}\in(0,\infty)$ such that
$$
\frac{c_{1}}{n}\leq\operatorname{ldel}_{n-1}\leq c_{2}.
$$
Theorem \ref{vertex} has been proved by choosing the vertices 
of a random polytope on the boundary of an Euclidean ball
whose radius is slightly bigger than $1$, namely, its radius is finely calibrated to be of the order
$1+\frac{1}{N^{\frac{2}{n-1}}}$.

A more general result, approximating a convex body by an arbitrary positioned polytope where the vertices are chosen with respect to a probability measure $P_f$ was proved in \cite{GW2018}. Approximation by arbitrary positioned polytopes in other metrics can be found in \cite{GTW2021,HSW2018}.

\section{Approximation of a polytope by a polytope with
fewer vertices}

B\'{a}r\'{a}ny conjectured in \cite{barany} that for every $n\in\mathbb N$
there is a constant $c_{n}\in(0,\infty)$ such that for every polytope
$P$  in $\R^n$ with a sufficiently large number of vertices $N$  there exists a vertex $x\in\R^n$ of $P$ such that the polytope $Q$, which is the convex hull of all
vertices of $P$ except $x$, satisfies
\begin{equation}
\label{eq1-1} \frac{\operatorname{vol}_n(P)-\operatorname{vol}_n(Q)}{\operatorname{vol}_n(P)} \leq c_{n}
N^{-\frac{n+1}{n-1}}\,.
\end{equation}
 B\'{a}r\'{a}ny also pointed out there that this
result implies a theorem of  Andrews \cite{andrews}, which
states that for a lattice polytope $P$ in $\R^{n}$ with $N$
vertices and positive volume, one has
$$
N^{\frac{n+1}{n-1}} \leq c_{n}
\operatorname{vol}_{n}(P).
$$
This conjecture has been confirmed by Reisner, Sch\"utt, and Werner in \cite{ReisSW}.

\begin{theorem}
\label{th3-2} \cite{ReisSW}
There are constants $c_0,c_2\in(0,\infty)$ such that
for every
$0<\epsilon<1/2$  the following holds: let $P$ be a polytope in
$\R^{n}$ having $N$ vertices $x_{1},\ldots,x_{N}$ with $N>
c_2^n/\epsilon$. Then there exists a subset $A\subset
\{1,\ldots,N\}$, with 
$\operatorname{card}(A)\geq (1-2\epsilon)N$, such that for
all $i\in A$: 
$$ 
\frac{\operatorname{vol}_{n}(P)
-\operatorname{vol}_{n}([x_{1},\ldots,x_{i-1},x_{i+1},\ldots,x_{N}])}
{\operatorname{vol}_{n}(P)} \leq
c_0\,n^{2}\epsilon^{-\frac{n+1}{n-1}}N^{-\frac{n+1}{n-1}}\,. 
$$
If $P$ is centrally symmetric, then we may replace $n^{2}$ by
$n^{3/2}$ on the right-hand side of the above inequality.
\end{theorem}

Let $x_{1},\dots,x_{N}$ be the vertices of $P$. 
In order to prove Theorem \ref{th3-2} we estimate the number of 
vertices $x_{i}$ whose distance from the polytope 
$[x_{1},\dots,x_{i-1},x_{i+1},\dots,x_{N}]$ is less than
$N^{-\frac{2}{n-1}}$.
We use inequality (\ref{ApproxHaus1}).

\section{Random approximation of polytopes by polytopes}

We turn now to question how well a random polytope of a polytope 
approximates this polytope.
B\'ar\'any and Buchta \cite{BB} studied how well a random polytope
whose vertices are chosen with respect to the uniform distribution on a polytope
approximates this polytope. It turns out that a random polytope $P_{N}$
of $N$ chosen points with respect to the uniform measure on a polytope $P$
satisfies
\begin{equation}\label{PolyRandom1}
\operatorname{vol}_n(P)-\E \operatorname{vol}_n(P_N) 
= \frac{ \operatorname{flag}(P)}{(n+1)^{n-1} (n-1)!} N^{-1} (\ln N)^{n-1} (1+o(1)),
\end{equation}
where $ \operatorname{flag}(P)$ is the number of flags of the polytope $P$.  The phenomenon that the expression should only depend on this combinatorial structure of the polytope, i.e., the flag number of the polytope, had been discovered
in connection with floating bodies by Sch\"utt (Theorem \ref{PolyRandomPoly}).
\par
Choosing random points from the interior of a convex body always produces a simplicial polytope with probability one. Yet often applications of the above mentioned results in computational geometry, the analysis of the average complexity of algorithms, and optimization necessarily deal with non-simplicial polytopes and it became crucial to have analogous results for random polytopes without this very specific combinatorial structure.
 \par
In all these results there is a general scheme: if the convex sets are smooth then the number of faces and the volume difference behave asymptotically like powers of $N$,  if the sets are convex polytopes then logarithmic terms show up. Metric and combinatorial quantities only differ by a factor $N$.
\vskip 3mm
Now we are discussing the case that the points are chosen from the boundary 
of a polytope $P$. This produces random polytopes which are neither simple nor simplicial and thus our results are a huge step in taking into account  the first point mentioned above. The applications in computational geometry, the analysis of the average complexity of algorithms, and optimization need formulae for the combinatorial structure of the involved random polytopes and thus the question on the number of facets and vertices are of interest. 
\par

To our big surprise  the volume difference contains no logarithmic factor. This is in sharp contrast to the results for random points inside convex sets.

\begin{theorem}\label{th:Vol}\cite{ReitzSW}
For the expected volume difference between a simple polytope $P \subset \R^n$ and the random polytope $P_N$ with vertices chosen uniformly at random from the boundary of $P$, we have
$$
\E (\operatorname{vol}_n(P)-\operatorname{vol}_n(P_N)) =
c_{n,P} N^{- \frac n{n-1}} (1+ O(N^{- \frac 1{(n-1)(n-2)}}))
$$ 
with some $c_{n,P}\in(0,\infty)$.
\end{theorem}

Intuitively, the difference volume for a random polytope whose vertices
are chosen from the boundary should be smaller than the one whose
vertices are chosen from the body. The previous result confirms this. The first one is of the order $N^{-\frac{n}{n-1}}$ compared to $N^{-1} (\ln N)^{n-1}$.
It is well known that for uniform random polytopes in the interior of a convex set the missed volume is minimized for the ball \cite{Bl2, Gr1, Gr2}, a smooth convex set, and -- in the planar case -- maximized by a triangle \cite{Bl2, DL, Gi} or more generally by polytopes \cite{BL}. Hence, one should also compare the result of Theorem \ref{th:Vol} to the result of choosing random points on the boundary of a smooth convex set. This  clearly leads to a random polytope with $N$ vertices. And by results of Sch\"utt and Werner \cite{SW5} and also Reitzner \cite{Re4}, the expected volume difference is of order 
$N^{- \frac 2{n-1}}$, which is smaller than the order in \eqref{PolyRandom1} as is to be expected, but also surprisingly much bigger than the order $N^{- \frac {n}{n-1}} $ occurring in Theorem \ref{th:Vol}.

We give a simple argument that shows that the volume difference
between the cube and a random polytope is at least of the order
$N^{-\frac{n}{n-1}}$. We consider the cube $C^{n}=[0,1]^{n}$
and the subset of the boundary
\begin{equation}\label{eq:simpl-cube}
\partial C^{n}\cap H^{+}\left((1,\dots,1),\left(\frac{(n-1)!}{n N}\right)^{\frac{1}{n-1}}\right)
=\bigcup_{i=1}^{n}\left(\frac{(n-1)!}{n N}\right)^{\frac{1}{n-1}}
[0,e_{1},\dots,e_{i-1},e_{i+1},\dots,e_{n}],
\end{equation}
which are the union of small simplices in the facets of the cube close to the vertices. Then 
$$
\frac{1}{N}=\l_{n-1} \left(\bigcup_{i=1}^{n}\left(\frac{(n-1)!}{n N}\right)^{\frac{1}{n-1}}
[0,e_{1},\dots,e_{i-1},e_{i+1},\dots,e_{n}]\right)
$$
and the probability
that none of the points $x_{1},\dots,x_{N}$ is chosen from this set
equals
$$
\left(1-\frac{1}{N}\right)^{N}\sim \frac{1}{e}.
$$
Therefore, with probability approximately $\frac{1}{e}$
the union of the simplices in \eqref{eq:simpl-cube}
is not contained in the random polytope and the difference volume is greater than
$$
\frac{1}{n!}\left(\frac{(n-1)!}{n N}\right)^{\frac{n}{n-1}}
\approx \frac{1}{n}N^{-\frac{n}{n-1}} ,
$$
which is in accordance with Theorem \ref{th:Vol}.

\subsection*{Acknowledgement}
JP is supported by the Austrian Science Fund (FWF) Project P32405 \textit{Asymptotic geometric analysis and applications} and Austrian Science Fund  (FWF)  Project  F5513-N26, which is a part of a Special
Research Program ``Quasi-Monte Carlo Methods: Theory and Applications''. EMW is partially supported by NSF grant DMS-2103482 and by a Simons Fellowship.

\bibliographystyle{plain}
\bibliography{bestandrandom}

\vskip 4mm

Joscha Prochno \\
  {\small        Universit\"at Passau}\\
    {\small      Fakult\"at f\"ur Informatik und Mathematik}\\
   {\small        Innstra{\ss}e 33}\\
    {\small      94032 Passau, Germany} \\
          {\small \tt joscha.prochno@uni-passau.de} \\

   Carsten Sch\"utt \\
     {\small       Christian Albrechts Universit\"at \hskip 20mm   Case Western Reserve University}\\
        {\small    Mathematisches Seminar \hskip 29mm   Department of Mathematics}\\
        {\small    24098 Kiel, Germany\hskip 35mmCleveland, Ohio 44106, U. S. A.} \\
         {\small \tt schuett@math.uni-kiel.de}   \\

     \and Elisabeth Werner\\
{\small Department of Mathematics \ \ \ \ \ \ \ \ \ \ \ \ \ \ \ \ \ \ \ Universit\'{e} de Lille 1}\\
{\small Case Western Reserve University \ \ \ \ \ \ \ \ \ \ \ \ \ UFR de Math\'{e}matiques }\\
{\small Cleveland, Ohio 44106, U. S. A. \ \ \ \ \ \ \ \ \ \ \ \ \ \ \ 59655 Villeneuve d'Ascq, France}\\
{\small \tt elisabeth.werner@case.edu}\\ \\

\end{document}